\documentclass[a4paper,12pt]{amsart}
\usepackage[english]{babel}
\usepackage[T1]{fontenc}
\usepackage{array}
\usepackage{makecell}
\usepackage[a4paper,margin=3cm]{geometry}
\usepackage[justification=centering]{caption}

\usepackage{amsmath,amsthm,amssymb,amsxtra,calligra,mathrsfs}
\usepackage{graphicx}
\usepackage{tikz}
\usepackage{tikz-cd} 
\usepackage{xcolor}
\usepackage{upgreek}
\usepackage{comment}
\usepackage{graphicx}
\usepackage{mathtools}
\usepackage{extarrows}
\usepackage{faktor}
\usepackage{enumitem}
\usepackage{comment}
\usepackage{bookmark}
\usepackage{hyperref}
\usepackage{mathrsfs}
\usepackage{multirow}

\hypersetup{
    colorlinks=true, 
    citecolor=blue,
    linkcolor=magenta,
    pdfauthor={},
	bookmarksnumbered=true,
}

\newcommand{\bF}{\mathbb{F}}

\newcommand{\bG}{\mathbb{G}}

\newcommand{\sheafHom}{\mathscr{H}\text{\kern -3pt {\calligra\large om}}\,}

\newcommand{\tabincell}[2]{\begin{tabular}{@{}#1@{}}#2\end{tabular}}

\DeclareMathOperator{\codim}{codim}

\DeclareMathOperator{\GL}{GL}
\DeclareMathOperator{\SL}{SL}
\DeclareMathOperator{\SO}{SO}
\DeclareMathOperator{\Sym}{Sym}

\DeclareMathOperator{\Gr}{Gr}

\DeclareMathOperator{\NL}{\mathrm{NL}}

\newcommand{\bZ}{{\mathbb Z}}

\newcommand{\bC}{{\mathbb C}}
\newcommand{\bQ}{{\mathbb Q}}
\newcommand{\bP}{{\mathbb P}}

\newcommand{\fX}{{\mathfrak X}}

\DeclareMathOperator{\scP}{\mathrm{P}}

\DeclareMathOperator{\seg}{seg}
\DeclareMathOperator{\pic}{Pic}
\DeclareMathOperator{\aut}{Aut}

\DeclareMathOperator{\ns}{NS}

\newcommand{\q}{/\!\!/\!}
\newcommand{\cO}{{\mathcal O}}
\newcommand{\cX}{{\mathcal X}}

\newcommand{\cY}{{\mathcal Y}}
\newcommand{\cI}{{\mathcal I}}

\newcommand\proj{\text{\rm Proj}}

\newcommand\im{\text{\rm Im}}

\newcommand{\cF}{{\mathcal F}}

\newcommand\cD{{\mathcal{D}}}
\newcommand\cS{{\mathcal{S}}}

\newcommand\cE{{\mathcal{E}}}

\newcommand\cL{\mathcal{L}}
\newcommand\cQ{{\mathcal{Q}}}
\newcommand\cC{{\mathcal{C}}}

\newtheorem{defn-pro}{Definition-Proposition}
\newtheorem{defn-thm}{Definition-Theorem}
\newtheorem{thm}{Theorem}[section]
\newtheorem{lem}[thm]{Lemma}
\newtheorem{cor}[thm]{Corollary}
\newtheorem{pro}[thm]{Proposition}

\newtheorem{defn}[thm]{Definition}

\newtheorem{question}[thm]{Question}

\newtheorem{conj}[thm]{Conjecture}
\newtheorem{rem}[thm]{Remark}

\theoremstyle{remark}

\title{Chow ring of moduli spaces of quasi-polarised K3 surfaces in lower genus}

\author{Fei Si}
\address{Beijing International Center for Mathematical Research, Peking University, No. 5 Yiheyuan Road Haidian District, Beijing, P.R.China 100871}
\email{sifei@bicmr.pku.edu.cn}

\begin{document}
\begin{abstract}
In the paper, we show the Chow ring of moduli spaces of quasi-polarised K3 surfaces in lower genus ($\le 5$) is tautological.  
\end{abstract}

\maketitle
\section{Introduction}
The intersection theory of moduli space $M_g$ of curves of genus $g$ has received much attention over the last decades,  especially after Mumford’s  work\cite{Mukai88} which introduced  the study of the tautological ring. This area remains an active
direction in algebraic geometry, as evidenced by the ICM survey \cite{Pandharipande18} and \cite{pixton21}.  A natural two dimensional generalisation  for  curves of genus $g$ is a quasi-polarised K3 surfaces of genus $g$, that is, a K3 surface with a big and nef primitive  line bundle $L$, where a general section of $L$ is a  curves of genus $g$. Similar to the case of curves
of genus $g$, the moduli stack  of quasi-polarised K3 surfaces of genus $g$ is a smooth  Delinge-Mumford stack and has coarse moduli space  $\cF_g$ with orbifold point at worst.
 The intersection theory on $\cF_g$ began with the work of the work of Marian-Oprea-Pandharipande \cite{MOP1}, where they first introduced the tautological ring $R^\ast(\cF_g)\subset A^\ast(\cF_g)$ as a $\bQ$-subalgebra of chow ring with $\bQ$-coefficient generated by pushforward of the kappa classes (see Definition \ref{defkappa}) from lower dimensional moduli spaces.  Later  \cite{BLMM17} Bergeron-Li-Millson-Moeglin's work showed that $R^1(\cF_g) \cong \pic(\cF_g)_\bQ $ and also give a beautiful formula for the dimension $\dim_\bQ R^1(\cF_g)$. Pandharipande-Yin \cite{PY20} proved that the tautological ring $R^\ast(\cF_g)$ is isomorphic to the subring $\NL^\ast(\cF_g)$ generated by Noether-Lefchetz locus. We refer to  \cite{Li20} for survey of intersection theory on moduli of K3 surfaces and hyperkahler varieties.

Recently,  Canning-Oprea-Pandharipande \cite[Theorem 1,(i)]{COP23} showed that $A^\ast(\cF_2)=R^\ast(\cF_2)$. Using the geometry of projective model  in genus $\le 5$, we generalise this result  for the moduli space of genus $g\le 5$.
\begin{thm}\label{chowtaut}
    The Chow ring  $A^\ast(\cF_g)$ is tautological, i.e, $A^\ast(\cF_g)=R^\ast(\cF_g)$  for $g \le 5$.  
\end{thm}

\begin{rem}
This can be viewed as an analogous result for the tautological property of the Chow ring of $M_g$ for small $g\leq 9$, s shown in the work of  \cite{Faber90a} \cite{Faber90b} \cite{Izadi95} \cite{FL08} \cite{PV15} and also recent work \cite{CL21} for $g=7,8,9$. 
\end{rem}

The idea is  to decompose the moduli space $\cF_g$ into $\cF_g=\cF_g^\circ \sqcup \mathop{\cup} \limits_{l\le g,d} \scP^g_{l,d}$ where $\cF_g^\circ $ is an open subset  and  $\mathop{\cup} \limits_{l\le g,d} \scP^g_{l,d}$ are union of some Noether-Lefschetz (NL) divisors, see Definition \ref{defNL}.  The geometric input in lower genus is that  the projective model for these K3 surfaces are clear. This enable us to  realise $\cF^\circ$ as a open subset in the quotient of the parameter space 
$P$ under a reductive group $G$ and the parameter space has the structure as a Grassmannian bundle. Via the localization sequence and the Grassmannian bundle constructions mentioned above, we known the generators of equivariant Chow groups for the parameter space. The locus we chose for the parameter space $P$ will imply the stabilizer of $G$ at each point in $P$ is finite and thus these generator for equivariant chow groups will decent to the generator for chows groups for $\cF_g^\circ$.  Then by the geometric construction of the parameter space  we  are able to show  these generators are tautological and then deduce the chow ring of $\cF_g^\circ $ is tautological. K3 surfaces parametrized by Noether-Lefschetz (NL) divisors $\scP^g_{l,d}$  also have nice projective model  and then the same trick will allow us to show the cycles on $\scP^g_{l,d}$  are also tautological. Applying localization sequence again, we obtain the results.

We hope the strategy works for moduli spaces of primitively quasi-polarised  K3 surfaces with Mukai models.  Recall the Mukai \cite{Mukai88} has classified a general K3 surface $(X,L) \in \cF_g$ for $g\le 12$ as a complete intersection in a homogeneous space or a Fano $3$-fold listed in the Table \ref{tab:Mukai}.  There is a natural GIT compactification $\overline{\cF}_g^{\rm GIT}$ for the open locus $\cF_g^\circ$.

\begin{table}[ht]
\renewcommand\arraystretch{1.5} 
\centering
\caption{Mukai Model for $(S,L)$}
\begin{tabular}{|c|c|c|c|} \hline 
$g$ &  General members in $\cF_g$ &  $\overline{\cF}_g^{\rm GIT}$ & $\cF_g-\cF_g^{\circ}$ \\\hline
3 &  quartic surface in $\bP^3$  & $|\cO_{\bP^3}(4)|\q \SL(4)$  & $\scP_{1,1}^3, \scP_{2,1}^3$ \\\hline
4 &   \tabincell{c}{cubic hypersurface on a smooth\\ quadratic threefold $Q$} & $|\cO_Q(3)|\q \SO(5)$   & $\scP^4_{1,1}, \scP^4_{2,1}, \scP^4_{3,1}$  \\\hline
5 &  \tabincell{c}{complete intersection of \\ three quadratics in $\bP^5$ } &  $\Gr(3, 21)\q \SL(6)$ &   $\scP^5_{1,1}, \scP^5_{2,1}, \scP^5_{3,1}$   \\ \hline 
6 &   \tabincell{c}{quadratic hypersurface on the \\ unique GM $3$-fold  } & $|\cO_Y(2)| \q \mathrm{PSL}(2)$  & \tabincell{c}{$\scP^{6}_{1,1}, \scP^{6}_{2,1}, \scP^{6}_{3,1}$ \\ $\scP_{4,1}^{6}, \scP^{6}_{5,2}$} \\\hline
7 &  \tabincell{c}{complete intersection of eight \\ hyperplanes in  the \\ isotropic Grassmannian  $\mathrm{IGr}(5,10)$}    & $\Gr(8, 16) \q  \mathrm{Spin}(10)$ & \tabincell{c}{$\scP^{7}_{1,1}, \scP^{7}_{2,1}, \scP^{12}_{3,1}$ \\ $\scP^{7}_{4,1}, \scP^{7}_{5,2},\scP^{7}_{6,2}$} \\\hline
8 &   \tabincell{c}{complete intersection of six \\ hyperplanes in $\Gr(2,6)$}  &  $\Gr(6, 15)\q \SL(6) $  & \tabincell{c}{$\scP^{8}_{1,1}, \scP^{8}_{2,1}, \scP^{8}_{3,1}$ \\ $\scP^{8}_{4,1}, \scP^{8}_{6,2},\scP^{8}_{7,2}$} \\\hline
9 &   \tabincell{c}{complete intersection of  four\\  hyperplanes  in Langrangian \\ Grassmannian $\mathrm{LGr}(3,6)$}  &  $\Gr(4,9)\q \mathrm{Sp}(6) $  & \tabincell{c}{$\scP^{9}_{1,1}, \scP^{9}_{2,1}, \scP^{9}_{3,1}$ \\ $\scP^{9}_{4,1}, \scP^{9}_{5,1},\scP^{9}_{6,2}$ \\  $\scP^{9}_{7,2}$  } \\\hline
10&  \tabincell{c}{complete intersection of three  \\ hyperplanes in flag variety associated \\with    the
adjoint representation of $G_2$}  &  $\Gr(3, 10)\q (G_2/Z(G_2))$ &   \tabincell{c}{$\scP^{10}_{1,1}, \scP^{10}_{2,1}, \scP^{10}_{3,1},\scP^{10}_{4,1}$ \\ $ \scP^{10}_{5,1},\scP^{10}_{7,2},\scP^{10}_{8,2}, \scP^{10}_{9,3}$  }\\\hline
12 & \tabincell{c}{ hyperplane section of a\\
smooth Fano threefold of degree 22}  & \tabincell{c}{$\bP^{13}$-\hbox{bundle over}\\ $\Gr(3,21 ) \q \SL(7)$ } & \tabincell{c}{$\scP^{12}_{1,1}, \scP^{12}_{2,1}, \scP^{12}_{3,1},\scP^{12}_{4,1},$ \\ $ \scP^{12}_{5,1},\scP^{12}_{6,1}, \scP^{12}_{7,2},\scP^{12}_{8,2}$ \\  $\scP^{12}_{9,3}$ $\scP^{12}_{10,3}$ } \\\hline
\end{tabular}  
\label{tab:Mukai}
\end{table}

The method also works for some lattice quasi-polarised K3 surfaces. In particular, 
we show 
\begin{thm}
   Let $N \hookrightarrow L_{K3}$ be a sublattice with gram matrix \[ \left( 
    \begin{array}{cc}
       0  & 2  \\
       2 & 0 
    \end{array}   
    \right) \ \ \hbox{or }\ \   \left( 
    \begin{array}{cc}
       2  & 2  \\
       2 & 0 
    \end{array}  
    \right)\]
    and $\cF_N$ the moduli space of $N$-quasi polarised K3 surfaces, then  chow ring $A^\ast(\cF_N)$ is tautological.
\end{thm}

\begin{conj}(Oprea-Pandharipand, see \cite{COP23}) \label{conjop} Let $N\hookrightarrow L_{K3}$ be lattice of signature $(1,r)$ and $\cF_N$ the moduli space of $N$-quasi-polarised K3 surfaces, then
\[   R^{19-r}(\cF_N) =0\ \hbox{and}\ \  R^{18-r}(\cF_N)=0 . \]
\end{conj}

The conjecture is widely open. The first evidence is due to Petersen. Let $RH^\ast(\cF_g)$ be the image of tautological ring under the cycle class maps $cl: A^\ast(\cF_g) \rightarrow H^\ast(\cF,\bQ)$. Then the cohomological version  $RH^{19-r}(\cF_N)=RH^{18-r}(\cF_N)=0$   for Conjecture \ref{conjop} is proved by Petersen \cite{Petersen19}.  For the conjectural vanishing in chow groups it is known in genus $g=2$ case by \cite{COP23} via explicit computations.  We remark that if the moduli space   of  lattice quasi-polarised K3 surfaces  is unirational,  
then $R^{19-r}(\cF_N) =0$ holds (see Remark \ref{rem:0cycle}). But the conjectural vanishing $R^{18-r}(\cF_N) =0$ seems much harder. Our result reduces the problem to show $A_1(\cF_g)=0$ for $g\le 5$.

\begin{rem}
    For $g >61$, $\cF_g$ is of general type by \cite{GHS07} and the chow ring $A^\ast(\cF_g)$ should be much more complicated.  It is expected $R^\ast(\cF_g) \subset A^\ast(\cF_g)$ is a proper subgroup. To find an explicit class in $ A^\ast(\cF_g)- R^\ast(\cF_g)$ is an interesting problem suggested by \cite{COP23}.
\end{rem}

\subsection*{Acknowledgement}
The author would like to thank Prof. Zhiyuan Li and Qizheng Yin for useful comments.

\subsection*{Structure of the paper}
In section \ref{set1}, we first review equivariant chow theory due to Totaro and Edidin-Graham.  Then we review the moduli space of (or lattice) quasi-polarised K3 surfaces and describe the projective geometry of K3 surfaces of genus $g\le 5$ and also two specific lattice quasi-polarised K3 surface.  In section \ref{set2}, we first show the chow rings  of  open part for the moduli spaces of K3 surfaces in lower genus are tautological.  This is based on the projective model geometry of K3 surfaces in section \ref{sec2.2}. Then we apply the same trick to show the NL divisor part and also some moduli spaces of lattice quasi-polarised K3 surfaces have  tautological chow ring. In the last section \ref{sec4} we discuss a possible geometric approach to deal with other $\cF_g$ with Mukai model, which is based on a stratification for  $\cF_g$ via a conjectural K-moduli wall-crossing result to connect GIT compactification $\overline{\cF}_g^{\rm GIT}$ in the Table \ref{tab:Mukai}  and Baily-Borel comactification for $\cF_g$.
\subsection*{Convention}
Throughout the paper, we will use the following convention.
\begin{itemize}
    \item $A^\ast(-)$   the Chow group with $\bQ$-coefficients.
    \item   $L_{K3}=\left( 
    \begin{array}{cc}
       0  & 1  \\
       1 & 0 
    \end{array}   
    \right)^3\oplus E_8^2$ the K3 lattice. 
 \end{itemize}

\section{Preliminaries} \label{set1}
\subsection{Equivaraint Chow ring }
We review the equivaraint chow  theory due to \cite{EG} \cite{Totaro}. Let 
$G$ be a reductive algebraic group over $k$. For any $i$, by \cite{Totaro}, there is a $n$-dimensional representation $V/k$  of $G$ such that there is a Zariski open subset $U\subset V$ where $G$ acts freely and $\codim (V-U) \geq d-i$. Assume $X$ is a scheme with $G$-action,  then
\begin{equation}
    X \times U \rightarrow X_G:= ( X \times U )/G
\end{equation}
is a principal $G$-bundle and $X_G$ exists as a scheme of finite type since $G$ is reductive. If $X$ is a point, $X_G$ can be viewed as an algebraic approximation of classifying space $BG$, by abuse of the notation, we also denote the algebraic approximation  $BG=U/G$, which is called Totora's approximation for $G$.

\begin{defn}
Notation as above, the $i$-th equivariant Chow groups for $d$-dimensional scheme $X$ is defined to be
\begin{equation}\label{eqch}
    \begin{split}
        A^G_i(X):=& A_{i+n-m}(X_G)=A_{i+n-m}(( X \times U )/G) \\
        =& A^{d-i}((X \times U )/G)
    \end{split}
\end{equation}
where $m=\dim G$. If $X$ is smooth ,  $  A_G^{d-i}(X):=A^G_i(X)$. 
\end{defn}
\begin{rem}
  $ A_G^{d-i}(X)$ can also be defined for any $X$ as equivariant operational chow group and it has ring structure. Moreover,  in smooth quotient case, it is isomorphic to $A^G_i(X)$. In this paper, we only deal with the smooth case $X$ and thus we adopt this definition for simplicity. 
\end{rem}

\begin{pro} (see \cite[\S 15 and\S 16]{Totaro}, also \cite{Pandharipande98}) \label{echow}
Let $G$ be an algebraic group $G\le \GL_n$ in the Table \ref{table:echow} and $V$ the standard representation of $\GL(n)$.  $V$ will induce a natural vector bundle over $B  \GL_n$ and we still denote the vector  bundle $V \rightarrow B\GL_n$. Let $c_i$ be $i$-th Chern class of vector bundle $V$, then the integral chow ring for classifying space is generated by the Chern classes $c_i$ of $V$ with possible relation given in the Table \ref{table:echow}.
\end{pro}
\begin{table}[ht]
    \centering
    \begin{tabular}{ | c | c|}
    \hline
Group $G$ & integral equivariant chow ring    \\ \hline
  $\GL_n$&  $\bZ[c_1,\cdots,c_n]$   \\ \hline
   $\SL_n$&  $\bZ[c_2,\cdots,c_n]$   \\ \hline
    $Sp_n$&  $\bZ[c_2,c_4\cdots,c_{2n}]$   \\ \hline
$O_n$&  $\bZ[c_1,\cdots,c_n]/(2c_1,2c_3,\cdots )$  \\ \hline
\end{tabular} 
  \caption{Equivariant Chow ring for some classical group }
  \label{table:echow}
\end{table}

\begin{rem}
    In the example listed above, indeed the cycle class map will induce isomorphism $A^\ast(BG) \rightarrow H^\ast(BG,\bQ)$. This is indeed true for any connected algebraic group.  For a general reductive group, it seems unknown.
\end{rem}

A very useful result for the chow ring of  nice quotient is the following and it will be frequently used in the next sections. 
\begin{pro}\label{finite}
Assume the stabilizer $G_x$ of each point $x\in G$  is finite group, then $A_G^\ast(X)=A^\ast(X/G)$.
\end{pro}

\subsection{K3 surface in lower genus and their moduli}\label{sec2.2}
A quasi-polarised K3 surface $(X,L)$ of genus $g$ consists a K3 surface $X$ with a big and nef line bundle $L$ such that $c_1(L)\in H^2(X,\bZ)$ is primitive and $L^2=2g-2>0$.  Let $N$ be an even lattice with signature $(1,r)$ and admitting an embedding $N \hookrightarrow L_{K3}:=\left( 
    \begin{array}{cc}
       0  & 1  \\
       1 & 0 
    \end{array}  
    \right)^{\oplus 3} \oplus E_8^2$. A more general notion is lattice $N$-quasi polarised K3 (see \cite[\S 1]{Dolgachev96}), which  means a K3 surface $X$ such that $N \hookrightarrow \ns(X)$ and there is a big and nef class $L \in N$.  Then  we can consider moduli space $\cF_N$ of lattice $N$-quasi-polarised K3 surfaces, i.e., $\cF_N$  parametrizes K3 surfaces $X$ admitting an embedding $N \hookrightarrow \pic(X)=\ns(X)$. By the famous Torelli theorem, $\cF_N= \widetilde{O}(N^\perp) \setminus \cD_{N^\perp}$ where $\cD_{N^\perp}$ is the period domain associated to the signature $(2,19-r)$ lattice $N^\perp$ and $\widetilde{O}(N^\perp)$ is the stable orthogonal group, i.e., the subgroup of orthogonal group $(N^\perp)$ whose elements act trivially on the discriminant group $A_{N^\perp}$.
    
\begin{defn}(NL divisors on $\cF_g$)\label{defNL}  A Noether-Lefthetz (NL)divisor on $\cF_g$ is the locus of K3 surfaces $(X,L)\in \cF_g$ with additional curve class $\beta$ such that 
\begin{equation}\label{NLlattice}
    \left( 
    \begin{array}{c|cc}
         &  L & \beta  \\ \hline  
       L     & 2g-2  & d \\
        \beta & d  & 2l-2
    \end{array}
    \right)
\end{equation}
is a primitve sublattice of $L_{K3}$.
\end{defn}
The NL divisor is irreducible by \cite[Theorem 2.1]{Litian21}   and we denote its reduced part by $\scP^g_{d,l}$.  As a sub locally symmetric space of codimension $1$, $\scP^g_{d,l}$ can be identified as the  quotient $\cD_{N^\perp} \cap v^\perp$ for a particular vector $v\in N^\perp$ under a finite index subgroup of $O((v)^\perp_{N\perp})$ but containing $\widetilde{O}((v)^\perp_{N\perp})$, in particular, there is a finite covering morphism 
\begin{equation}
    \cF_N \rightarrow \scP^g_{d,l}.
\end{equation}
where the rank two lattice has gram matrix as in  (\ref{NLlattice}).

A key input for the understanding the chow ring of moduli space of quasi-polarised K3 of genus $g\le 5$ is  the  projective model of these surfaces. Now let us describe them for $g\le 5$.
\subsubsection{$g=3$ case} \label{decom3}
By   the work of Saint-Donat \cite{Saint-Donat74}, a general K3 surface $(X,L) \in \cF_3$
is a quartic surface $S\subset \bP^3$. Denote $\cF_3^{\circ} \subset \cF_3$ the locus of quartic surface with simple singularities. Then  $\cF_3-\cF_3^{\circ}$ consists of two NL divisor. By the work of \cite{LO21} and \cite{ADL22}
connecting the GIT  space $|\cO_{\bP^3}(4)|\q \SL_4$ of quartic surface and Baily-Borel compactification $\cF_3^\ast$,  the following  result is obtained 
\begin{pro}\label{decomp3}
 There is a decomposition  $\cF_3=\cF_3^{\circ} \sqcup \scP_{2,1} \sqcup \scP_{1,1}$ for $\cF_3$  such that 
 \begin{enumerate}
     \item $(X,L) \in \cF_3^{\circ}$ if and only if a quartic surface $S\subset \bP^3$ with simple singularities.
     \item $(X,L) \in \scP_{2,1}$ if and only if $X$ is  a double cover $X \rightarrow \bP^1 \times \bP^1$ branched along bidegree $(4,4)$ curve. In this case, $\scP_{2,1}$ is  called hyperelliptic divisor. 
     \item $(X,L) \in \scP_{1,1}$ if and only if $X$ is an ellipic K3 surface $X \rightarrow \bP^1$ such that $X\subset \bP (\cO_{\bP^1} \oplus \cO_{\bP^1}(8) \oplus \cO_{\bP^1}(12))$. In this case, $\scP_{1,1}$ is  called unigonal divisor. 
 \end{enumerate}
\end{pro}
 \begin{proof}
    (1) is obvious. For (2),  by the K-moduli space  wall-crossing interpolation  of  GIT space $|\cO_{\bP^3}(4)|\q \SL_4$ and Baily-Borel compactification $\cF_3^\ast$, $(X,L)\in \scP_{2,1}$ if and only if $X$ is a double cover of del pezzo surface  $X$ of degree $8$ such that $X$ is degeneration of $\bP^1\times \bP^1$ and $(X,(\frac{1}{2}-\epsilon)C)$ is a K-polystable del pezzo pair where $ C\sim -2K_X$. By \cite{ADL23}, for all such $(X,(\frac{1}{2}-\epsilon)C)$ del Pezzo pairs, underlying surface $X$ is $\bP^1 \times \bP^1$ and the curve $C \subset \bP^1 \times \bP^1$ is bi-degree $(4,4)$ curve. For (3), it is a direct consequence that a K3 surface in $\scP_{1,1}$  admits a elliptic fibration with a section and has the Weierstrass model.
 \end{proof}   

\subsubsection{$g=4$ and $g=5$ case} \label{decom45}
Similar to the K3 surface of genus $3$ case, due to the work of Saint-Donat
a general K3 surface $(X,L)$ of genus $4$ or $5$ is a complete intersection of quadric and cubic in $X\subset \bP^4$ or complete intersection of three quadrics  $X \subset \bP^5$,  see also \cite[\S 2]{Litian21}. For $g=4$ or $5$,  denote $\cF_g^\circ$  the locus of K3 surface of genus $g$, which is complete intersections $\bP^g$.  A  crucial geometric result is that \cite[Lemma 2.2]{Litian21}  or  see also \cite[Proposition 7.15 and Proposition  7.19]{Saint-Donat74}.
\begin{lem} \label{tri45}
   Let $g=4$ or $5$, then  a   K3 surface $(X,L) \in \scP_{3,1}^g$  if and only if  
   \begin{enumerate}
       \item   the line bundle $L$ will induce a birational map  $X \dashrightarrow \bP^4$ with image a  $(2,3)$ complete intersection $X \subset \bP^4$ where the quadric is singular ;
       \item   for $g=5$, 
    \begin{enumerate}
        \item  either  the line bundle $L$ will induce a birational map  $X \dashrightarrow \bP^5$ with image a bidegree $(2,3)$  hypersurface in $\bP^1 \times \bP^2$ 
        \item or there is decomposition for curve class $L=L_0+E+R$  where $L_0$ is a curve of genus $2$, $C$ is a rational curve and $E$ is an irreducible elliptic curve with 
\begin{equation}\label{trimatrix}
     \left( 
    \begin{array}{c|cc c}
         &  L_0 & E  & R  \\ \hline  
       L_0     & 2  &  2 & 1 \\
        E& 2  &    0  &  1\\
        R & 1  & 1 & -2
    \end{array}
    \right) .
\end{equation}
    \end{enumerate}
   \end{enumerate}
\end{lem}

\begin{pro} \label{decomp45} Notation as above, 
    \begin{enumerate}
     \item  for $g=4$, there is a decomposition $\cF_g=\cF_g^\circ \sqcup \scP_{2,1} \cup \scP_{1,1}$  where 
     \begin{enumerate}
         \item $(X,L) \in \scP_{2,1}$ in and only if $X$ is a double cover $S \rightarrow \bF_i$ branched along a curve $C\sim -2K_{\bF_i}$ for $i=1$ or $4$. In this case, $(X,L)$ is called hyperelliptic K3 and lies in NL divisor $\scP_{2,1}$
         \item $(X,L) \in \scP_{1,1}$ in and only if $X$ is an ellipic K3 surface with section.
     \end{enumerate}

     \item  for $g=5$, there is a decomposition $\cF_g=\cF_g^\circ \sqcup \scP_{3,1} \cup  \scP_{2,1} \cup \scP_{1,1}$  where  
     \begin{enumerate}
         \item  $(X,L) \in \scP_{3,1}$ if and only  if either  $X$ is a bidegree $(2,3)$  hypersurface in $\bP^1 \times \bP^2$ or $X \rightarrow \bP^2$ is double over  branched along a nodal sextic curve, which is tangent to a line.
       \item   a double cover $S \rightarrow \bP^1 \times \bP^1$ branched along bidegree $(4,4)$ curve. In this case, $(X,L)$ is called hyperelliptic K3 and lies in NL divisor $\scP_{2,1}$
       \item $(X,L) \in \scP_{1,1}$ in and only if $X$ is an ellipic K3 surface with section.
     \end{enumerate}
       \end{enumerate}  
\end{pro}
\begin{proof}
For $g=4$, by Lemma \ref{tri45} (1),  the complement of $\cF_g^\circ$ consists of two NL divisors $\scP_{2,1}$ and $\scP_{1,1}$. If  $(X,L) \in \scP_{1,1}$, then $X$ is an ellipic K3 surface with section.   Assume $(X,L) \in \scP_{2,1}$, by the K-moduli  wall-crossing result in \cite[Theorem 1.5]{PSW23} connecting Baily-Borel compactification of $\scP_{2,1}$ and GIT space for sections in $|-2K_{\bF_1}|$, $\cF_N$ is a big  open subset of K-moduli space of del Pezzo pairs of degree $8$, where all singular del Pezzo pairs are classified. In particular,  $\cF_N$  be decomposed into $\cF_N^\circ\sqcup  H_u$ where 
\begin{itemize}
    \item $\cF_N^\circ$ parametrizes  K3 surfaces  $X\rightarrow \bP^2$ branched along a sextic curve with at least a node point and  $X$ has ADE singularities at worst,
    \item $H_u$ parametrizes K3 surface $X \rightarrow \bF_4$ branched along $C \sim -2K_{\bF_4}$ so that $X$ has ADE singularities at worst.
\end{itemize}

For $g=5$, by Lemma \ref{tri45} (2), if $(X,L) \in \scP_{3,1}$ is not s a bidegree $(2,3)$  hypersurface in $\bP^1 \times \bP^2$ , then  its Neron-Severi group contains a sublattice with the following gram matrix  
\[\left( 
    \begin{array}{c|cc c}
         &  L_0 & L_0-E  & R  \\ \hline  
       L_0     & 2  &  0 & 1 \\
         L_0-E& 0  &    -2  &  0\\
        R & 1  & 0 & -2
    \end{array}
    \right). \]
Then $L_0$ will induce a degree $2$ morphism $X \rightarrow \bP^2=\bP H^0(X,L_0)^\vee$. It is a double cover branched along $B \sim -2K_{\bP^2}$ such that the rational curve $ L_0-E$ contracts to $A_1$-singular point on $B$, while rational $R$ maps to a line tangent to the branched curve $B$.  This proves part (a) of (2). For (b), just note that there is a finite onto morphism  $ \cF_N \rightarrow  \scP^{g}_{2,1}$ for $N= \left( 
    \begin{array}{cc}
       0  & 2  \\
       2 & 0 
    \end{array}   
    \right) $ and such lattice $N$-quasi-polarised K3 surfaces are described in the proof of Proposition \ref{decomp3}.
\end{proof}
\subsection{Tautological rings and Noether-Lefchetz locus }
We review the tautological ring $R^\ast(\cF_g)$ defined in  \cite{MOP1}. Similar definition for moduli space of hyperkahler varieties is given by  \cite{BL19}.
\begin{defn}\label{defkappa}
Let   $\pi: (\fX,\cL_1,\cL_2,\cdots,\cL_r) \rightarrow B$ be the universal family of quasi-polarised K3 surfaces over the base $B$.  The kappa class is defined     \[\kappa_{a,b}(\pi):= \pi_\ast(  c_1(\cL_1)^{a_1}\cdot c_1(\cL_2)^{a_2}\cdot c_1(\cL_r)^{a_r}\cdot c_2(T_{\pi})^b )\  \in\  A^{a_1+\cdots+a_r+2b-2}(B)\]
\end{defn}

Assume $N$  is the lattice of signature $(1,r)$ admits a class $\ell \in N$ with $\ell^2=2g-2$ and $N\hookrightarrow L_{K3}$ then it induces a primitive embedding  $N^\perp \hookrightarrow \ell^\perp \hookrightarrow L_{K3}$ and it turns out there is morphism $\iota_N: \cF_N \rightarrow \cF_g$ which is finite onto image, known as Nother-Lefthetz locus on $\cF_g$.   Denote $\NL^\ast(\cF_g) \subset A^\ast(\cF_g)$ the $\bQ$-sub algebra generated by Nother-Lefthetz locus, that is, the  pushfoward  of fundamental class of $\cF_N$ under $\iota_N: \cF_N \rightarrow \cF_g$.  The tautological ring $R^\ast(\cF_g) \subset A^\ast(\cF_g)$ is defined as  the $\bQ$-sub algebra generated by the pushforward of all kappa classes $\iota_{N \ast} (\kappa_{a,b})$ on $\cF_N$.  Clearly, $\NL^\ast(\cF_g) \subset R^\ast(\cF_g)$. 
A fundamental theorem proved by Pandharipande-Yin is the following
\begin{thm}(\cite{PY20})
    $\NL^\ast(\cF_g) = R^\ast(\cF_g)$
\end{thm}
For the $\bQ$-Picard group of $\cF_g$, Bergeron-Li-Millson-Moeglin’s work \cite{BLMM17} shows all divisor classes are tautological 
\begin{thm}\label{bletc}
  $\pic(\cF_g)_\bQ$ is generated by Noether-Lefchetz divisors.  In particular, $R^1(\cF_g)= A^1(\cF_g) \cong \pic(\cF_g)_\bQ$.
\end{thm}

\begin{rem}
The above construction of  tautological rings can be parallel generalised to the moduli spaces of quasi-polarised hyperkahler varieties. Then  a cohomological version for the above result is proved by Bergeron-Li \cite{BL19} using the method of automorphic forms.  
\end{rem}

\section{Chow ring of lattice quasi-polarised K3 surfaces }
\label{set2}
\subsection{Chow ring of open part} \label{sec3.1}
In this subsection, we are going to show 
\begin{thm}
The Chow ring $A^\ast(\cF_g^\circ)$ is generated by  tautological classes for $g\le 5$.
\end{thm}
The proof will be divided into the following discussion.
\subsubsection{$g=3$ case} Let $V=\bC^4$ be the standard representation for $\SL_4$ and $E:= \Sym^4(V^\vee)$.
By Proposition  \ref{decomp3}, 
$\cF_g^\circ$ is the coarse moduli space $U/G$ of the quotient stack $[U/G]$ where $U \subset \bP E$ is the parameter space of quartic surfaces with ADE singularities. Since each point in $U$ has finite stabilizer under the action $G$ since the quotient stack $[U/G]$ is a DM stack \footnote{see \cite[Chapter 4]{Huybrechts16}  for the argument to show the quotient stack parametrizing ADE K3 surfaces is a DM stack.}    or by shah's GIT stability characterization \cite{Shah81} in this case.  Then by the result from  equivariant chow theory in  Proposition \ref{finite}, \[ A^\ast(\cF_g^\circ)=A^\ast(U/G)=A^\ast_G(U).\]
Applying localization sequence to the open embedding $U \hookrightarrow \bP E$, we get 
\begin{equation*}
    A_G^\ast(Z) \xrightarrow{i_{G,\ast}} A^\ast_G(\bP E) \rightarrow A^\ast_G(U)=A^\ast(\cF_g^\circ) \rightarrow 0
\end{equation*}
where $Z:=\bP E-U$  is the complement of $U$. After base change to Totaro's approximation space $BG$,  there is also morphism \begin{equation*}
    p: (\bP E)_G \rightarrow BG
\end{equation*} and it is a projective bundle over the $BG$. By the construction, $A^\ast_G(\bP E)=A^\ast ( (\bP E)_G)$. 
Thus by the bundle structure, we have
\begin{equation}
    A^\ast_G(\bP E)=A^\ast ( (\bP E)_G)= \frac{p^\ast A^\ast(BG)[H]}{I_G}
\end{equation}
where $H:=c_1(\cO_{(\bP E)_G}(1))$ is the tautological divisor class and \[ I_G= 
\big \langle \mathop{\sum} \limits_{i=0}^{r} c_i(p^\ast E) H^{r-i} \big \rangle  \] is the Grothendieck relation. Since $E=\Sym ^4(V^\vee)$,   $c_i(p^\ast E)$ is a homogeneous polynomial in $c_1(p^\ast V^\vee),\cdots c_4(p^\ast V^\vee) $ of degree $i$.  By Proposition \ref{echow}, $p^\ast A^\ast(BG)$ can be identified with the Chern class $c_i(p^\ast V^\vee)$.  

In a summary, the above discussion shows $A^\ast(\cF_g^\circ)$ is generated by class  $H$ and $c_1(p^\ast V^\vee),\cdots c_4(p^\ast V^\vee) $.  

\begin{pro}
 The classes $c_i(p^\ast V^\vee) $ restricts to $A_G^\ast(U)=A^\ast(U/G)$  are tautological.    The class $H|_{U}$ is propositional to Hodge line bundle on $\cF_g^\circ$, in particular, it is tautological class.
\end{pro}
\begin{proof}
Let  $j: \fX \hookrightarrow  \bP p^\ast V$ be the universal hypersurface over $(\bP E)_G$ and $f: (\fX,\cL) \rightarrow (\bP E)_G$  the flat projective family by projection onto $(\bP E)_G$  where $\cL:=\cO_{\bP p^\ast V}(1)|_{\fX}$ is $f$-ample. Then \[ \cX \sim \pi^\ast H+c_1(\cO_{\bP p^\ast V}(4))\    \]  in $A^1(\bP p^\ast V)$ and there is diagram  
\begin{equation*}
    \begin{tikzcd}
        \fX \arrow[d,"f"] \arrow[r,hook,"j"] & \bP p^\ast V^\vee \arrow[dl,"\pi"]  \\
       (\bP E)_G  & 
    \end{tikzcd}
\end{equation*}
 Note that $\cL$ is relatively very ample and    $f_\ast \cL \cong p^\ast V^\vee $  up  to twist some power line bundle $H$ on $(\bP E)_G$.
  By the Grothendieck-Riemann-Roch formula,
  \begin{equation*} \begin{split} ch(p^\ast V^\vee)=ch(\pi_{\ast} \cL)=ch(\pi_{!} \cL ) = \pi_\ast ( \mathop{\sum} \limits_{n} \frac{c_1(\cL)^n}{n!} \cdot td(T_\pi)) \end{split}\end{equation*}
    where right-hand side is linear combination of kappa class and left-hand side is the polynomial in  $c_1(p^\ast V^\vee),\cdots c_4(p^\ast V^\vee) $. 
This proves  $c_1(p^\ast V^\vee),\cdots c_4(p^\ast V^\vee) $ is tautological. Note that $U/G \subset \bP E \q G$ and the Picard number of $\cF_g^\circ =U/G$ is $1$.  Thus the decent of $H|_U$ must be a positive  multiple of Hodge line bundle on $\cF_g^\circ$ since Hodge line bundle is ample. It is well known  
 Hodge line bundle is tautological. 
\end{proof}

\begin{rem}
    By the solution Maulik-Pandharipande's conjecture  due to  Bergeron-Li-Millson-Moeglin using arithmetic methods,   each divisor class in $A^1(\cF_g) \cong \pic(\cF_g)_\bQ $ is tautological. 
\end{rem}
\subsubsection{$g=4$ case}

Let $V$ be a standard representation of $G=\SL_5$. By the geometric input for $\cF_g^\circ$ in  \S \ref{decom45}, there is a parameter space $\bP E \rightarrow  \bP \Sym^2(V^\vee)$ for $(2,3)$ complete section in $\bP^4=\bP V$. To apply equivariant chow theory,  we give the construction over  $BG$, which mimics usual construction.   Let $g: \bP \Sym^2(V^\vee) \rightarrow BG$ be projective bundle over $BG$. Then there is a family  $Q\subset \bP g^\ast V^\vee $  of quadric surface over $\bP \Sym^2(V^\vee)$, which is a general section in $|\cO_{\bP \Sym^2(V^\vee)}(1)\otimes \cO_{ \bP g^\ast V^\vee }(2)|$  and denote $q: Q  \rightarrow \Sym^2(V^\vee)$. Then the pushforward $\cE:=  q_{\ast}\cO_{\bP g^\ast V^\vee}(3)$ is a vector bundle of rank  $30$  on $\bP \Sym^2(V^\vee)$. The projective bundle $\bP \cE \rightarrow \bP \Sym^2(V^\vee)$ is  a tower of bundle over Totaro's approximation space $BG$. By the construction, $\bP \cE=(\bP E)_G$.  As before, we get intermediately 
\begin{equation*}
  A^\ast_G(\bP E)=A^\ast(\bP \cE)= \frac{p^\ast A^\ast(BG)[h,\xi]}{I}
\end{equation*}
where 
\begin{itemize}
    \item  $h:=c_1(\cO_{\bP \cE}(1))$ and $\xi:= c_1(\cO_{\bP \Sym^2(V^\vee)}(1))$ 
    \item  $I$ is an ideal given by  two Grothendieck  relation in $h$ and $\xi$ respectively.
\end{itemize}
\begin{pro}
The chow ring $A^\ast(\cF_g^\circ)$ is tautological for $g=4$.     
\end{pro}
\begin{proof}
   As before, it is sufficient to show the generators are tautological.  It is similar to deal with $(2,3)$ complete intersection in $\bP^3$ as in proof of Proposition \ref{hypertau}  and we omit the details here.
\end{proof}

\subsubsection{$g=5$ case}
Let $V$ be a standard representation of $G=\SL_6$. There is a natural $\SL_6$-action on  $P:=\Gr(3,  \Sym ^2(V^\vee))$. As before base change to $BG$,  then there is a Grassmannian bundle 
\[  p: P_G \rightarrow BG \]
where $ P_G $ is a Grassmannian bundle over $BG$.
There is universal subbundle and quotient bundle of $p^\ast \Sym ^2(V^\vee)$ in the following 
\begin{equation}\label{universalbd}
    0 \rightarrow \cS \rightarrow p^\ast \Sym ^2(V^\vee) \otimes \cO \rightarrow \cQ \rightarrow 0 .
\end{equation}
By the basis theorem for Grassmannian bundle in \cite[Proposition 14.6.5]{Ful}, \[A^\ast(P_G) \cong p^\ast A^\ast(BG)[ c_i(\cS), c_j(\cQ)] / I\] where $I$ is the relation given by 
\begin{equation*}
     c(\cS) c(\cQ)=p^\ast c(\Sym ^2(V^\vee)).
\end{equation*}
In particular, localization sequence argument will show $A^\ast(\cF_g^\circ)$  is generated by the restriction of class $c_2(p^\ast V^\vee),\cdots c_5(p^\ast V^\vee) $ and $c_{i}(\cS )$.
There is a universal family of $\fX \subset \bP p^\ast  V^\vee $ with the projection $f: \fX \rightarrow P_G$ with 
 diagram 
 \begin{equation*}
     \begin{tikzcd}
         \fX \arrow[d,"f"]  \arrow[r,hook] &  \bP p^\ast V^\vee  \arrow[dl] \\
          P_G&
     \end{tikzcd}
 \end{equation*}
 There is short exact sequence 
\begin{equation} \label{universal}
     0 \rightarrow \cI_\fX \rightarrow \cO \rightarrow \cO_{\fX} \rightarrow 0.
\end{equation}
By flatness of $f$ and direct computation of stalk, we have $R^1f_\ast (\cI_\fX \otimes  \cO_{ \bP p^\ast V^\vee}(2)) =0$.
 Thus by twisting line bundle $\cO_{ \bP p^\ast V^\vee}(2)$ with (\ref{universal}) and then taking the pushfoward, we will give exact sequence on $P_G$ 
\begin{equation}\label{push1}
    \begin{split}
        0 &\rightarrow f_\ast (\cI_\fX \otimes \cO_{ \bP p^\ast V^\vee}(2)) \rightarrow  f_\ast ( \cO_{ \bP p^\ast V^\vee}(2))  \rightarrow  f_\ast (\cO_{\fX} \otimes \cO_{ \bP p^\ast V^\vee}(2)) \rightarrow 0
    \end{split}
\end{equation}
  Now we claim  the  exact sequence  in (\ref{push1}) coincide with the universal bundle (\ref{universalbd}). Clearly, $f_\ast ( \cO_{ \bP p^\ast V^\vee}(2))  \cong \Sym^2( p^\ast V^\vee)$. Note that the geometric interpretation of surjective bundle map $\Sym^2( p^\ast V^\vee) \rightarrow f_\ast (\cO_{\fX} \otimes \cO_{ \bP p^\ast V^\vee}(2))$ is evaluation map on each fiber  and thus the kernel sheaf is $f_\ast (\cI_\fX \otimes \cO_{ \bP p^\ast V^\vee}(2))$.
\begin{pro}
The chow ring $A^\ast(\cF_g^\circ)$ is tautological for $g=5$.   
\end{pro}
\begin{proof}
First we show $c_i(p^\ast V^\vee)$. This follows  Grothendieck-Riemann-Roch formula for the relative very ample line bundle  $\cO_{\bP p^\ast V^\vee}(1)|_{\fX}$ directly.  Then we deal with the Chern classes of $\cS$. A key point is that by the geometric construction,  as a cycle on $\bP p^\ast V^\vee$, 
\begin{equation}\label{familyrel}
    \fX \sim  ( c_1(\cO_{\bP p^\ast V^\vee}(2))+c_1(\cS))^3 
\end{equation} 
    The Grothendieck-Riemann-Roch formula  combining vanishing $Rf^{>0}_\ast (\cI_\fX \otimes \cO_{ \bP p^\ast V^\vee}(2))=0$ shows 
    \begin{equation*}
        \begin{split}
             ch(\cS)&=ch(f_\ast (\cI_\fX \otimes \cO_{ \bP p^\ast V^\vee}(2)))=ch(f_\ast (\cI_\fX \otimes \cO_{ \bP p^\ast V^\vee}(2) td(T_f))\\
             &=f_\ast (ch(\cO_{\fX}) td(T_f)).
        \end{split}
    \end{equation*}
 Then the class relation (\ref{familyrel}) implies  the Chern classes of universal subbundle $\cS$ are tautological.   
\end{proof}
\subsection{Chow ring of NL divisors $\scP_{d,1}^g$.}
We follow the similar strategy to show the Chow ring of NL divisors in \S \ref{decom3} and \ref{decom45} are tautological. 
\begin{thm}
 The Chow ring $A^\ast(\scP^g_{d,1})$ is tautological for $g\le 5$ and $d \le 3$.   
\end{thm}

\subsubsection{Unigonal divisor }
Note that the unigonal divisor $\scP^g_{1,1}$ always 
 admit a surjective morphism  $\cF_N \rightarrow \scP^g_{1,1}$ from moduli space $\cF_N$  of  special lattice  polarised K3 whose gram matrix is 
 \[  N=\left( 
    \begin{array}{cc}
       -2  & 1  \\
       1     &  0
    \end{array}
    \right) ,\]
which is exactly  the moduli space of ellipic K3 with a section. Such K3 surface admits Weierstrass model, that is, it is a hypersurface in ruled $3$-fold $\bP (\cO_{\bP^1} \oplus \cO_{\bP^1}(8) \oplus \cO_{\bP^1}(12)) $  which is cut out by Weierstrass equation 
\[ y^2z=x^3+axz^2+bz^3\]
where $(a,b) \in  H^0(\cO_{\bP^1}(8))\oplus H^0(\cO_{\bP^1}(12)$. Via this geometric input, the moduli space of ellipic K3 surface can be identified as quotient space $ W/\SL_2 \times \bG_m$ where $W \subset H^0(\cO_{\bP^1}(8))\oplus H^0(\cO_{\bP^1}(12)$ is the locus of elliptic K3 with ADE singularities at worst.
By the work of \cite{CK23} and see also \cite[\S 3.2]{COP23},  
\begin{pro}
  The chow ring for moduli space of  ellipic K3 surfaces is tautological. 
\end{pro}
As a direct corollary, 
\begin{cor}
$A^\ast(\scP^g_{1,1})$ is tautological.
\end{cor}

\subsubsection{ Chow ring of hyperelliptic divisor}
There are two types of hyperelliptic divisor $\scP^g_{2,1}$ given by moduli $\cF_N$ of $N$-quasi-polarised K3 surface with gram matrix 
 \[   \left( 
    \begin{array}{cc}
       0  & 2  \\
       2 & 0 
    \end{array}   
    \right) \ \ \hbox{and}\ \   \left( 
    \begin{array}{cc}
       2  & 2  \\
       2 & 0 
    \end{array}  
    \right) \ \  . \]
If $N= \left( 
    \begin{array}{cc}
       0  & 2  \\
       2 & 0 
    \end{array}   
    \right) $, a general K3 surface  in $\cF_N$  is 
a double cover of $\bP^1 \times \bP^1$ branched along bi-degree $(4,4)$-curves. The moduli space of bi-degree $(4,4)$ curves can be compactified by VGIT space of $(2,4)$ complete intersection: let $\cY \subset \bP^3 \times |\cO_{\bP^3 }(2)|$ be the universal quadratic hypersurface and $\pi_1: \cY \rightarrow  |\cO_{\bP^3 }(2)|$ and  $\pi_2: \cY \rightarrow  \bP^3$ the two projections. Then it is not hard to show $\cE:=\pi_{1 \ast} (\pi_2^\ast \cO_{\bP^3 }(4))$ is a vector bundle of rank $25$. Let $\pi: \bP \cE \rightarrow |\cO_{\bP^3 }(2)|=\bP^{9}$ be the projective bundle and  denote $h:=\pi^\ast c_1(\cO_{\bP^{9}}(1))$ and  $\xi=c_1(\cO_{\bP \cE }(1))$, then the VGIT space is defined as
\[ \bP \cE \q _{L_t} \SL(4) := \proj( R(\bP \cE ,L_t)^G  ),\ \   L_t=h+t \xi,\ \ 0< t<\frac{1}{2}.\] By the result of VGIT wall crossing in  \cite{LO21} or K-moduli wall crossing in \cite{ADL23},  the VGIT space will interpolate  GIT space $|\cO_{\bP^1 \times \bP^1}(4,4)|\q \SL(2) \times \SL(2)$  and  Baily-Borel's compactification for $\scP_{2,1}$. In particular,  $\scP^g_{2,1}$ can be identified as a quotient $W/\SL_4$ of a open subset $W$ in $\bP \cE$ under the action $\SL_4$. Here $W$ parametrizes $(2,4)$ complete intersection  where the double cover produces a K3 with ADE singularities at worst.

\begin{pro}\label{hypertau}
    The chow ring $A^\ast(\cF_N)$ is tautological for $N= \left( 
    \begin{array}{cc}
       2  & 2  \\
       2 & 0 
    \end{array}   
    \right) $.
\end{pro}
\begin{proof}
 Following the strategy as before , we  need to construct $G=\SL_4$-equivariant version for the parameter space  $(\bP \cE)_G \rightarrow BG$ of $(2,3)$ complete intersection. Let $V=\bC^4$ be the standard representation of $G=\SL_4$, $G=\SL_4$-equivariant version is just base change to the Totao's approximation space $BG$
 \begin{equation*}
     (\bP \cE)_G \xrightarrow{\pi} (\bP  \Sym^2(V^\vee))_G \rightarrow BG
 \end{equation*}
 Then \[A^\ast((\bP \cE)_G)= p^\ast A^\ast(BG)[h, \xi]/I\] where abuse of notation denote $h:=\pi^\ast c_1(\cO_{(\bP  \Sym^2(V^\vee))_G}(1))$ and $\xi:= c_1(\cO_{(\bP \cE)_G} (1))$.  There are universal family $(\fX,\cC)\subset \bP p^\ast V^\vee$ over $(\bP \cE)_G$ such that as a cycle on  $\bP p^\ast V^\vee$,  $\fX \sim  h + \cO_{\bP p^\ast V^\vee}(3)$  and $\cC\sim $.    The double covering $\varphi: \cY \rightarrow \fX$ branched along universal $(2,4)$ complete intersection curve $\cC$ is the universal hyperellitpci K3 surface over $(\bP \cE)_G$, these families are flat over $(\bP \cE)_G $ with the following commutative diagram 
 \begin{equation*}
     \begin{tikzcd}
    \cY  \arrow[dr,"f_1"] \arrow[rr,"\varphi"]&  & \fX \arrow[dl,""]   \arrow[r,hook,"j"]& \bP p^\ast V^\vee  \arrow[dll,"f_3"] \\
       & (\bP \cE)_G & &
     \end{tikzcd}
 \end{equation*}
Observe that $p^\ast A^\ast(BG)$ is the polynomial ring in \[p^\ast \seg_i(V^\vee)= \seg_i(p^\ast V^\vee)=f_{3\ast}(\cO_{\bP p^\ast V^\vee}(1)^{i+3}), \ i=2,3,4 \]
Since  $ \varphi^\ast  j^\ast \cO_{\bP p^\ast V^\vee}(1)$ is a $f_1$-ample, which gives a polarised family  of K3 surfaces $(\cY,\varphi^\ast  j^\ast \cO_{\bP p^\ast V^\vee}(1))$ over $ (\bP \cE)_G$, then  by the diagram and projection formula
\begin{equation}\label{push}
     \begin{split}
         f_{1\ast }(\varphi^\ast  j^\ast \cO_{\bP p^\ast V^\vee}(1))^m)=&  f_{3\ast } j_\ast \varphi_\ast (\varphi^\ast  j^\ast \cO_{\bP p^\ast V^\vee}(1))^m)\\
         =&2 h\cdot f_{3\ast } (\cO_{\bP p^\ast V^\vee}(1))^m)+6 f_{3\ast } (\cO_{\bP p^\ast V^\vee}(1))^{m+1}) \\
         =& 2 h\cdot \seg_{m-3}( p^\ast V^\vee)+6 \seg_{m-2}( p^\ast V^\vee).
     \end{split}
\end{equation}
The left-hand side in (\ref{push}) is the kappa class thus tautological. Since the divisor class $h$ is tautological  by Theorem \ref{bletc}, the induction on $m$ will imply all Segree classes $\seg_i(p^\ast V)$ are tautological. Thus the Chern classes  are tautological since $\seg_i(V^\vee)=(-1)^i c_i(V^\vee)$.
\end{proof}

If $N= \left( 
    \begin{array}{cc}
       0  & 2  \\
       2 & 0 
    \end{array}   
    \right) $, the geometry of K3 surface in $\cF_N$
are similar to the case  $N= \left( 
    \begin{array}{cc}
       2  & 2  \\
       2 & 0 
    \end{array}   
    \right) $. Each K3 surfaces in $\cF_N$ is described as in Proposition \ref{decom45} (1) (a).
\begin{pro}
    The chow ring $A^\ast(\cF_N)$ is tautological.
\end{pro}
\begin{proof}
 The  proof is also  based on the  geometric construction of the hyperelliptic family  and  the same argument as in \S \ref{sec3.1}. We leave it to the interested readers. 
\end{proof}

\subsubsection{Tri-ellitpic divisor}
Let $ W\subset |\cO_{\bP^1 \times \bP^2}(2,3)|$ be the locus of bidegree $(2,3)$ hypersurfaces with ADE singularities and $G=\SL_2 \times \SL_3$ acting on $W$. Let $V_2$ (resp. $V_3$) be standard representation of $\SL_2$ (resp. $\SL_3$). Denote $\bP \Sym^2(V_2^\vee) \otimes \Sym^3(V_3^\vee)$ the projective bundle over $BG$, then $W_G \subset \bP \Sym^2(V_2^\vee) \otimes \Sym^3(V_3^\vee)$ is $W \times_G BG$ and thus the localization sequence implies $A^\ast(W/G)=A^\ast(W_G)$ is generated by class in $p^\ast A^\ast(BG)$  and $H$ where $p: \bP \Sym^2(V_2^\vee) \otimes \Sym^3(V_3^\vee)$ and $H:=c_1(\cO_{\bP \Sym^2(V_2^\vee) \otimes \Sym^3(V_3^\vee)}(1))$ since  \[A^\ast(\bP \Sym^2(V_2^\vee) \otimes \Sym^3(V_3^\vee))=\frac{p^\ast A^\ast(BG)[H]}{I}.\]
\begin{pro} \label{triopen}
    The Chow ring $A^\ast(W/G)$ is tautological. 
\end{pro}
\begin{proof}
 It is sufficient to show the classes in $p^\ast A^\ast(BG)$ are tautological.  One thing need to note that in our case $G=\SL_2 \times \SL_3$, by \cite[Proposition 14.2]{Totaro}, $A^\ast(BG)$ is generated by class $p^\ast c_i(V_2)$ , $p^\ast c_i(V_3)$.    Thus it is also sufficient to show $p^\ast c_i(V_2)$ , $p^\ast c_i(V_3)$ and $H$ are tautological.  One can also use the universal family $\fX$ over $\bP \Sym^2(V_2^\vee) \otimes \Sym^3(V_3^\vee)$, but the construction is little different from previous case: let $ \bP p^\ast V_i^\vee$ be projective bundle over $\bP \Sym^2(V_2^\vee) \otimes \Sym^3(V_3^\vee)$  for $i=2,3$ and then take the fiber product  \[ \bP:= \bP p^\ast V_2 \times_{\bP \Sym^2(V_2^\vee) \otimes \Sym^3(V_3^\vee)} \bP p^\ast V_3 \]  of the two projective bundles. There is a flat projective family $\fX \subset \bP$ of K3 surfaces $\pi: \fX \rightarrow \bP \Sym^2(V_2^\vee) \otimes \Sym^3(V_3^\vee)$   with the following diagram 
  \begin{equation*}
      \begin{tikzcd}
        \bP \arrow[rr,""]  \arrow[dd,""]   &  & \arrow[dd,""]  \bP p^\ast V_3 \\
        &   \fX \arrow[dr,"\pi"] \arrow[dl,"\mu_1"]\arrow[ur,"\mu_2"] \arrow[ul,hook]& \\
          \bP p^\ast V_2  \arrow[rr,""]&    & \bP \Sym^2(V_2^\vee) \otimes \Sym^3(V_3^\vee).
      \end{tikzcd}
  \end{equation*}
  Then for $a,b\in \bZ_{\ge 0}$ there is natural morphism of coherent sheaves
 \begin{equation}\label{bundlemor}
      \pi_\ast(\mu_1^\ast \cO_{\bP p^\ast V_2}(a)) \otimes  \pi_\ast(\mu_2^\ast \cO_{\bP p^\ast V_3}(b)) \rightarrow   \pi_\ast(\mu_1^\ast \cO_{\bP p^\ast V_2}(a) \otimes \mu_2^\ast \cO_{\bP p^\ast V_3}(b)). 
 \end{equation}
  Since both $R\pi_\ast^{>0}(\mu_1^\ast \cO_{\bP p^\ast V_2}(a))=0$ and $R\pi_\ast^{>0}(\mu_1^\ast \cO_{\bP p^\ast V_3}(b))=0$ by the construction,  the morphism (\ref{bundlemor}) is isomorphic after stalkwise checking.  Then applying Grothendieck-Riemann-Roch formula 
  \begin{equation*}
     \begin{split}
         &ch(\pi_\ast(\mu_1^\ast \cO_{\bP p^\ast V_2}(a) \otimes \mu_2^\ast \cO_{\bP p^\ast V_3}(b)))\\
         =&ch(\pi_!(\mu_1^\ast \cO_{\bP p^\ast V_2}(a) \otimes \mu_2^\ast \cO_{\bP p^\ast V_3}(b)))\\
         =& \pi_\ast (ch(\mu_1^\ast \cO_{\bP p^\ast V_2}(a) ) ch(\mu_2^\ast \cO_{\bP p^\ast V_3}(b)) td(T_\pi)),
     \end{split}
  \end{equation*}
  thus chern classes of vector bundles $\pi_\ast(\mu_1^\ast \cO_{\bP p^\ast V_2}(a)=\Sym^a(p^\ast V_2)$ and $\pi_\ast(\mu_2^\ast \cO_{\bP p^\ast V_2}(a)=\Sym^a(p^\ast V_3)$ are tautological. This finishes proof. 
\end{proof}

\begin{rem}
Here the NL divisor can be dominant by moduli space $\cF_N$ of lattice $ N=\left( 
    \begin{array}{cc}
       2 & 3  \\
       3 & 0 
    \end{array}  
    \right)$ quasi-polarised K3 surfaces. But some $N$-quasi-polarised K3 surface may have more complicated projective model. That is why we deal with $\scP_{3,1}$ directly, unlike the unigonal and hyperelliptic case.  A GIT model for $\cF_N$ is constructed in \cite{CWY22}. It would be interesting to find a resolution of birational period map from GIT to $\cF_N^\ast$ via K-moduli wall crossing,  and then  using the geometric interpolation to show whether $A^\ast(\cF_N)\cong R^\ast(\cF_N)$. See the strategy in \S \ref{sec4} for more details. 
\end{rem}
Now we deal with the remaining locus of $\scP_{3,1}^g$ described in part (a) of Proposition \ref{decom45} (2). Recall such K3 surface is a double cover of $\bP^2_{[x,y,z]}$ branched along a nodal sextic curve tangent to a line in $\bP^2$.  Under the action of $\SL_3$, we may assume the nodal point is $p=[0,0,1]$  and the line $l=\{z =0\}$, then we have the equation for such sextic curve of the form 
\[ z^4f_2(x,y)+z^3f_3(x,y)+z^2f_4(x,y)+zf_5(x,y)+g_2(x,y)=0\]
where $f_i,g_i$  are homogeneous polynomial in $x,y$ of degree $i$. The stabilizer group of the point $p$ and line $l$ is $G=\SL_2$ acting on $x,y$.  Denote $V$ a standard representation of $G=\SL_2$ and \[E:=\Sym^2(V^\vee) \oplus \Sym^3(V^\vee) \oplus \Sym^4(V^\vee) \oplus \Sym^5(V^\vee) \oplus \Sym^2(V^\vee)\] by the above discussion,  there is $G$-invariant  open locus $U \subset \bP E$ parametrizing sextic curves such that its double cover K3 surface has only ADE singularities. Then $A^\ast(U/G)=A^\ast(U_G)$ and the argument of equivariant chow theory as in previous cases will show $A^\ast(U/G)$ is tautological. By localization combining 
Proposition \ref{triopen}, we will have
\begin{cor}
    $A^\ast(\scP_{2,1}^g)$ is tautological for $g=5$. 
\end{cor}

  \vspace{0.3cm}

\subsubsection*{Proof of Theorem \ref{chowtaut}}  Recall $\cF_g^\circ \subset \cF_g$ is  the open locus whose complement $\cF_g-\cF_g^\circ$ is the union of NL divisors. Then there is a localization sequence 
\begin{equation*}
    A_\ast (\cF_g-\cF_g^\circ) \xrightarrow{i_\ast}  A_\ast(\cF_g) \rightarrow  A_\ast(\cF_g^\circ) \rightarrow 0
\end{equation*}
where \[i:\cF_g-\cF_g^\circ= \mathop{\cup} \limits_{l} \scP^g_{l,1} \hookrightarrow \cF_g\] is the closed embedding. The image of $i_\ast$ is also in the union of the image of pushforward $A_\ast(\scP^g_{l,1}) \rightarrow A_\ast(\cF_g)$ induced by the closed embedding $\scP^g_{l,1} \hookrightarrow \cF_g$.  Thus $A_\ast(\cF_g)$ is tautological if $A_\ast(U)$ is tautological and  $A^\ast(P_{l,1})$ is tautological for each irreducible component $P_{l,1}$ in $\cF_g-\cF_g^\circ$. 

\begin{rem}\label{rem:0cycle}
Let $N\hookrightarrow L_{K3}$ be an even lattice of signature $(1,r)$ and $\cF_N$ the moduli space of $N$ quasi-polarised K3 surfaces. If $\cF_N$ is unirational, then $R^{19-r}(\cF_N)=R_0(\cF_N)=A_0(\cF_N)=0$. For example, $\cF_N=\cF_g$  for $g \le 12$ or $N$ is taken as in \cite{BHK16}.  The argument for the vanishing $=A_0(\cF_N)=0$ is the following: 
Let $\cF^\ast_N$ be the Baily-Borel's compactification for $\cF_N$, then the boundaries $\cF^\ast-\cF_g$ consists of modular curves $C_1,\dots ,C_m$.  Note that unirationality of  $\cF^\ast_N$ implies that $\cF^\ast_N$  is rational connected.  It is known that  a  rational connected projective varieties  has  $A_0(\cF_g^\ast)=\bQ$. Then  by the localization sequence 
\[   \mathop{\bigoplus} A_0(C_i) \rightarrow   A_0(\cF_g^\ast ) \rightarrow   A_0(\cF_g) \rightarrow 0,\]
$A_0(\cF_g) =0$ follows as if the image of pushforward of $0$-cycles on $C_i$ is not identically Zero. This is clear since there is degree of $0$-cycles on $C_i$.
\end{rem}

\section{Further discussion}\label{sec4}

Let $\overline{P}^K_g(c)$ be the K-moduli space parametrizing S-equivalence class of $c$-K-semistable Fano $3$-fold pairs $(X,S)$ where $K_X+S \sim 0$ and $(-K_X)^3=r^3(2g-2)$ \footnote{Here we take the main irreducible component.}. Classification results for smooth Fano $3$-folds implies $g\le 10$ or $g=12$. By many people's work, $\overline{P}^K_g(c)$ exists as a projective variety. Moreover,   the work \cite{ADL19} shows there is wall-crossing behavior for the space $\overline{P}^K_g(c)$ when $c\in (0,1) \cap \bQ$ varies.  Observe that there is a natural forgetful rational map $\overline{P}^K_g(c) \dashrightarrow \cF_g$. Motivated by the HKL program for the moduli of surfaces,   we expect the following wall-crossing behavior holds 
\begin{enumerate}
    \item $\overline{P}^K_g(c) \cong \overline{P}_g^{GIT}$,  
    \item there is a dominant morphism $\overline{P}^K_g(1-\epsilon) \rightarrow \cF_g^\ast$. 
\end{enumerate}
This is confirmed in the case of   quartic surfaces pairs by \cite{ADL23} and also partially for quadric $3$-fold pairs in \cite{FLLST}.  A key observation is that (2) will ensure that there is a decomposition $\cF_g=\cF_g^\circ \sqcup (\cup \scP_{h,l})$ where $\cF_g^\circ \subset \overline{P}_g^{GIT}$ is listed in the Table \ref{tab:Mukai}, while K3 surface in  $\scP_{h,l}$ can be embedded into a specific Fano $3$-fold $Y$.  Then one can follow the similar strategy to show $A^\ast(\cF_g)$ is tautological for $g$ in the Table \ref{tab:Mukai}.  But one need to deal with the following question.
Let $Y$ be a Gorestein $3$-fold with $|-K_Y|$ empty and $G:=\aut(Y)$ reductive, denote 
 $W\subset |-K_Y|$ the locus of anti-canonical sections with ADE singularities, then there is $G$-equivariant family of $\cS \subset W \times Y$ of K3 surfaces over $W$. The equivariant chow theory construction as before produce a family K3 surface $\cS \rightarrow W_G$ such that $\cS\subset W_G \times Y$ and $W_G$ is viewed as variety  $p: W_G \rightarrow BG$ over Totaro's approximation space $BG$.
\begin{question}
    Is the image $\im (p^\ast: A^\ast(BG) \rightarrow A^\ast((W)_G)=A^\ast(W/G))$  tautological ?
\end{question}
Prior we do not know the generators of the equivariant chow ring $A^\ast(BG)$ for a reductive group $G$. If the answer is no, it will give the construction of non-tautological classes.

\vspace{1cm}
\bibliographystyle{alpha}
\bibliography{lowgenus}
\end{document}